\documentclass[lualatex]{article}
\usepackage{authblk}

\title{Convergence rates analysis of a multiobjective proximal gradient method}

\author[1]{Hiroki Tanabe}
\author[1]{Ellen H. Fukuda}
\author[1]{Nobuo Yamashita}

\affil[1]{Graduate School of Informatics, Kyoto University}
\affil[ ]{\{tanabehiroki@amp.i.kyoto-u.ac.jp\},\{ellen,nobuo\}@i.kyoto-u.ac.jp}

\date{}

\usepackage{graphicx}
\usepackage{amssymb}
\usepackage{mathtools}
\usepackage{epstopdf}
\usepackage{setspace}
\usepackage{amsmath}
\usepackage{enumitem}

\usepackage{amsthm}
\usepackage{cases}
\usepackage{xparse}
\usepackage{hyperref}
\usepackage[capitalize]{cleveref}

\def\usebfsetcapital{\def\setcapital##1{\mathbf{##1}}}

\usebfsetcapital

\def\setR{\setcapital{R}}

\newcommand\condition[1]{\quad \text{#1}}
\newcommand\forallcondition[1]{\condition{for all~$#1$}}

\DeclareMathOperator*{\argmin}{argmin}

\DeclareMathOperator{\dom}{dom}
\DeclareMathOperator{\prox}{prox}

\DeclarePairedDelimiter{\norm}{\lVert}{\rVert}

\DeclarePairedDelimiter{\set}{\lbrace}{\rbrace}
\DeclarePairedDelimiterX{\Set}[2]{\lbrace}{\rbrace}{#1\mathrel{}\delimsize\vert\mathrel{}#2}
\DeclarePairedDelimiterX{\innerp}[2]{\langle}{\rangle}{#1, #2}
\newcommand{\T}{\top\hspace{-1pt}}

{}

\newcommand{\tom}{\ensuremath{\set{1, \dots, m}}}

\newtheorem{theorem}{Theorem}[section]
\newtheorem{proposition}{Proposition}[section]
\newtheorem{lemma}{Lemma}[section]
\newtheorem{definition}{Definition}[section]

\newtheorem{remark}{Remark}[section]
\newtheorem{assumption}{Assumption}[section]
\newtheorem{algorithm}{Algorithm}[section]

\crefname{equation}{}{}
\Crefname{equation}{Eq.}{Eqs.}
\crefname{figure}{Figure}{Figures}
\crefname{assumption}{Assumption}{Assumptions}

\begin{document}
\maketitle
\begin{abstract}
    Many descent algorithms for multiobjective optimization have been developed in the last two decades. Tanabe et al. (Comput Optim Appl 72(2):339--361, 2019) proposed a proximal gradient method for multiobjective optimization, which can solve multiobjective problems, whose objective function is the sum of a continuously differentiable function and a closed, proper, and convex one. Under reasonable assumptions, it is known that the accumulation points of the sequences generated by this method are Pareto stationary. However, the convergence rates were not established in that paper. Here, we show global convergence rates for the multiobjective proximal gradient method, matching what is known in scalar optimization. More specifically, by using merit functions to measure the complexity, we present the convergence rates for non-convex ($O(\sqrt{1 / k})$), convex ($O(1 / k)$), and strongly convex ($O(r^k)$ for some~$r \in (0, 1)$) problems. We also extend the so-called Polyak-{\L}ojasiewicz (PL) inequality for multiobjective optimization and establish the linear convergence rate for multiobjective problems that satisfy such inequalities ($O(r^k)$ for some~$r \in (0, 1)$).
\end{abstract}

\section{Introduction} \label{sec: intro}
Let us consider the following unconstrained multiobjective optimization problem:
\begin{equation} \label{eq: MOP}
    \min_{x \in \setR^n} \quad F(x)
,\end{equation}
where~$F \colon \setR^n \to (\setR \cup \set*{\infty})^m$ is a vector-valued function with~$F \coloneqq (F_1, \dots, F_m)^\T$.
We assume that each component~$F_i \colon \setR^n \to \setR \cup \set*{\infty}$ is defined by
\begin{equation} \label{eq: f + g}
F_i(x) \coloneqq f_i(x) + g_i(x) \forallcondition{i \in \tom}
,\end{equation}
where~$f_i \colon \setR^n \to \setR$ is continuously differentiable and~$g_i \colon \setR^n \to \setR \cup \{\infty\}$ is closed, proper, and convex.
Now, assume that each~$\nabla f_i$ is Lipschitz continuous.

Many descent methods have been developed for multiobjective optimization problems~\cite{Fukuda2014}.
The objective value decreases at each iteration with these methods, and the convergence of sequences generated by them is established under reasonable assumptions.
For example, the steepest descent method~\cite{Fliege2000} converges globally to Pareto stationary points for differentiable multiobjective optimization problems.
Other multiobjective descent methods include projected gradient~\cite{Grana-Drummond2004}, proximal point~\cite{Bonnel2005}, Newton's~\cite{Fliege2009}, trust-region~\cite{Carrizo2016}, and conjugate gradient methods~\cite{LucambioPerez2018}.
Methods for infinite-dimensional vector optimization have also been studied, such as the proximal point~\cite{Bonnel2005} and the inertial forward-backward methods~\cite{Bot2018}. 
For the multiobjective problem~\cref{eq: MOP}, we can apply the proximal gradient method proposed in~\cite{Tanabe2019}.
Note that the steepest descent and the projected gradient methods can be considered as special cases of the proximal gradient method.
It was shown that each accumulation points of the sequences generated by the proximal gradient method are Pareto stationary under some assumptions, but their convergence rate has not been studied yet.
In this paper, we analyze the convergence rates of the multiobjective proximal gradient method by generalizing the related studies of single-objective gradient methods.

In recent years, there have been many studies about convergence rates of gradient methods for scalar optimization problems given by~\cref{eq: MOP} with~$m = 1$ and~$g_1 = 0$.
By using gradient methods, it is known that under mild assumptions $\left\{ \norm*{ \nabla f_1(x^k) } \right\}$ converges to zero with rate~$O(\sqrt{1/k})$ for non-convex cases~\cite{Beck2017,Nesterov2004}, $\{ f_1(x^k) - f_1^\ast \}$ converges to zero with rate~$O(1 / k)$ for convex problems, where~$f_1^\ast$ is the optimal objective value~\cite{Beck2017,Beck2009}, and $\{ x^k \}$ converges linearly to the optimal solution for problems with strongly convex objective functions~\cite{Beck2017,Nesterov2004}.
It is also shown that $\{ f_1(x^k) - f_1^\ast \}$ converges linearly to zero for problems that satisfy the so-called Polyak-{\L}ojasiewicz (PL) inequality~\cite{Karimi2016,Polyak1963}.
Moreover, similar results are shown for the proximal gradient method, which generalizes the gradient method, for a scalar optimization problem given by~\cref{eq: MOP} with~$m = 1$~\cite{Beck2017,Fukushima1981,Karimi2016}.
On the other hand, there are few studies about convergence rates for multiobjective gradient methods.
In~\cite{Fliege2019}, the authors analyzed the complexity of multiobjective steepest descent methods by showing that for convex cases $\{ \sum_i \lambda_i^k (F_i(x^k) - F_i(x^\ast)) \}$ converges to zero with rate~$O(1/k)$ for some constant~$\lambda_i^k$.
However, the constant $\lambda_i^k$ depends on the algorithm, so we cannot compare the convergence rate with other methods.
Reference~\cite{Bello-Cruz2020} estimated the number of iterations of the proximal gradient method for convex multiobjective optimization, but the issue related to the algorithm-dependent measure still remains.
Furthermore, they only analyze problems with convex objective functions and do not with non-convex or strongly convex objective functions.

In this paper, we analyze convergence rates of the multiobjective proximal gradient method by using merit functions~\cite{Tanabe2022a} to measure the complexity, which permits the comparison with other methods.
We also show a linear convergence rate for particular problems.
Note that we can apply these results to the multiobjective steepest descent method, too, since the multiobjective proximal gradient method is a generalization of them.

The outline of this paper is as follows.
We present some notations and notions of Pareto optimality and Pareto stationarity in \cref{sec: preliminaries}.
\Cref{sec: merit} recalls two merit functions for multiobjective optimization proposed in~\cite{Tanabe2022a}, and \Cref{sec: PL} generalizes PL inequality to the multiobjective optimization problem~\cref{eq: MOP}.
We estimate the global convergence rates of the proximal gradient method for multiobjective optimization in \cref{sec: convergence rates}, and finally, we conclude this paper in \cref{sec: conclusion}.

\section{Preliminaries} \label{sec: preliminaries}
\subsection{Notations}
We first present some notations that will be used in this paper.
Let us denote by~$\setR$ the set of real numbers.
We use the symbol~$\norm*{ \cdot }$ for the Euclidean norm in~$\setR^n$.
The notation~$u \le v$ ($u < v$) means that~$u_i \le v_i$ ($u_i < v_i$) for all~$i \in \tom$.
For problem~\cref{eq: MOP}, we also define the effective domain of~$F$ by $\dom F \coloneqq \{ x \in \setR^n \mid F(x) < \infty \}$.
In addition, let $F(A) \coloneqq \{F(x) \in \setR^m \mid x \in A\}$ and $F^{-1}(B) \coloneqq \{x \in \setR^n \mid F(x) \in B \}$ for $A \subseteq \setR^n$ and $B \subseteq \setR^m$.
Moreover, we call
\[
    h^\prime (x; d) \coloneqq \lim_{t \searrow 0} \frac{h(x + t d) - h(x)}{t}
\]
the directional derivative of~$h \colon \setR^n \to \setR \cup \{ \infty \}$ at~$x$ in the direction~$d$.
Note that~$h^\prime(x; d) = \nabla h(x)^\T d$ when~$h$ is differentiable at~$x$, where~$\nabla h(x)$ stands for the gradient of~$h$ at~$x$ and~$\T$ denotes transpose.
Furthermore, for a convex function~$\hat{h} \colon \setR^n \to \setR \cup \{\infty\}$, a vector~$\eta \in \setR^n$ is called subgradient of~$\hat{h}$ at a point~$x \in \setR^n$ if
\[
    \hat{h}(y) \ge \hat{h}(x) + \eta^\T (y - x) \forallcondition{y \in \setR^n}.
\]
The set of all subgradients of~$\hat{h}$ at~$x$ is denoted by~$\partial \hat{h}(x)$ and is called the subdifferential of~$\hat{h}$ at~$x$.

\subsection{Pareto optimality}
Now, we introduce the concept of optimality for the multiobjective optimization problem~\cref{eq: MOP}.
Recall that~$x^\ast \in \setR^n$ is \emph{Pareto optimal} if there is no~$x \in \setR^n$ such that~$F(x) \le F(x^\ast)$ and~$F(x) \neq F(x^\ast)$.
Likewise,~$x^\ast \in \setR^n$ is \emph{weakly Pareto optimal} if there does not exist~$x \in \setR^n$ such that~$F(x) < F(x^\ast)$.
It is known that Pareto optimal points are always weakly Pareto optimal, and the converse is not always true.
We also say that~$\bar{x} \in \setR^n$ is Pareto stationary~\cite{Custodio2011}, if and only if,
\[
    \max_{i \in \tom} F_i^\prime(\bar{x}; d) \ge 0 \forallcondition{d \in \setR^n}.
\]
We state below the relation among the three concepts of Pareto optimality.
\begin{lemma}~\cite[Lemma~2.2]{Tanabe2019}
    For problem~\cref{eq: MOP}, the following statements hold.
    \begin{enumerate}
        \item If~$x \in \setR^n$ is weakly Pareto optimal for~\cref{eq: MOP}, then~$x$ is Pareto stationary.
        \item Let every component~$F_i$ of~$F$ be convex.
              If~$x \in \setR^n$ is Pareto stationary for~\cref{eq: MOP}, then~$x$ is weakly Pareto optimal.
        \item Let every component~$F_i$ of~$F$ be strictly convex.
              If~$x \in \setR^n$ is Pareto stationary for~\cref{eq: MOP}, then~$x$ is Pareto optimal.
    \end{enumerate}
\end{lemma}

\section{Merit functions for multiobjective optimization} \label{sec: merit}
We introduce two merit functions for the multiobjective optimization problem~\cref{eq: MOP} proposed in~\cite{Tanabe2022a} to keep the paper self-contained.
Recall that each~$F_i \colon \setR^n \to \setR \cup \{\infty\}$ is the sum of a continuously differentiable function~$f_i$ and a closed, proper, and convex function~$g_i$, and that~$\nabla f_i$ is Lipschitz continuous with Lipschitz constant~$L_i > 0$.
Let us define 
\begin{equation} \label{eq: L}
    L \coloneqq \max_{i \in \tom} L_i
.\end{equation}
The first merit function is the simple function~$u_0 \colon \setR^n \to \setR \cup \{ \infty \}$ defined as follows:
\begin{equation} \label{eq: u_0}
    u_0(x) \coloneqq \sup_{y \in \setR^n} \min_{i \in \tom} \{ F_i(x) - F_i(y) \}
.\end{equation}
The second type is the regularized and partially linearized merit function~$w_\ell \colon \setR^n \to \setR$, given by
\begin{equation} \label{eq: w_ell}
    w_\ell(x) \coloneqq \max_{y \in \setR^n} \min_{i \in \tom} \left\{ \nabla f_i(x)^\T (x - y) + g_i(x) - g_i(y) - \frac{\ell}{2} \norm*{ x - y }^2 \right\}
,\end{equation}
where~$\ell > 0$ is a given constant.
\begin{remark} \label{rem: merit}
    When $m = 1$, we have
    \[
        \begin{dcases}
            u_0(x) = F_1(x) - F_1^\ast, \\
            w_\ell(x) = \frac{1}{2 \ell}\mathcal{D}_{g_1}(x, \ell) \ge \frac{\ell}{2} \norm*{x - \prox_{g_1/\ell} \left(x - \frac{1}{\ell} \nabla f_1(x) \right) }^2
,
        \end{dcases}
    \]
    where $F_1^\ast$ is the optimal objective value and
    \begin{equation} \label{eq: D_g}
        \begin{dcases}
        \mathcal{D}_{g_1}(x, \alpha) \coloneqq - 2 \alpha \min_{y \in \setR^n} \left\{ \nabla f_1(x)^\T (y - x) + \frac{\alpha}{2} \norm*{y - x}^2 + g_1(y) - g_1(x) \right\}, \\
        \prox_{h}(z) \coloneqq \argmin_{y \in \setR^n} \left\{ h(y) + \frac{1}{2} \norm*{y - z}^2 \right\}.
        \end{dcases}
    \end{equation}
    Note that when $g_1 = 0$, we have $\mathcal{D}_{g_1}(x, \alpha) = \norm*{\nabla f_1(x)}^2$.
    The functions~$u_0$ and $w_\ell$ are clearly merit functions for scalar-valued optimization, at least in terms of stationarity~\cite{Beck2017}.
\end{remark}
The following theorem shows that in multiobjective cases, $u_0$ and $w_\ell$ are also merit functions in the Pareto sense.
\begin{theorem}~\cite[Theorems 3.1 and 3.9]{Tanabe2022a} \label{thm: merit Pareto}
    Let~$u_0$ and~$w_\ell$ be defined by~\cref{eq: u_0} and~\cref{eq: w_ell}, respectively for all~$\ell > 0$.
    Then, the following statements hold.
    \begin{enumerate}
        \item For all $x \in \setR^n$, we have $u_0(x) \ge 0$.
            Moreover, $x \in \setR^n$ is weakly Pareto optimal for~\cref{eq: MOP} if and only if $u_0(x) = 0$.
        \item For all $x \in \setR^n$, we have $w_\ell(x) \ge 0$.
            Moreover, $x \in \setR^n$ is Pareto stationary for~\cref{eq: MOP} if and only if $w_\ell(x) = 0$.
    \end{enumerate}
\end{theorem}
From Theorem~\ref{thm: merit Pareto} and the lower semicontinuity of~$u_0$ and~$w_\ell$,~$u_0(x^k) \to 0$ and~$w_\ell(x^k) \to 0$ for a bounded sequence~$\set*{x^k}$ imply weak Pareto optimality and Pareto stationarity of the accumulation points, respectively.
In this paper, we evaluate the theoretical performance of the multiobjective proximal gradient method by considering convergence rates of~$\set*{u_0(x^k)}$ and~$\set*{w_\ell(x^k)}$, while literature~\cite{Tanabe2019} shows only Pareto stationarity of the accumulation points.
The following theorem connects the two merit functions~$u_0$ and~$w_\ell$.
\begin{theorem}~\cite[Theorems 4.1 and 4.2]{Tanabe2022a} \label{thm: merit relation}
    Let~$u_0$ and~$w_\ell$ be defined by~\cref{eq: u_0} and~\cref{eq: w_ell}, respectively, for all~$\ell > 0$.
    Then, the following statements hold.
    \begin{enumerate}
        \item Suppose that~$\nabla f_i$ is Lipschitz continuous with Lipschitz constant~$L_i > 0$, and let $L \coloneqq \max_{i \in \tom} L_i$. Then, we have
            \[
                u_0(x) \ge w_L(x) \forallcondition{x \in \setR^n}.
              \] \label{enum: u_0 >= w_L}
        \item Assume that each~$f_i$ is strongly convex with modulus~$\mu_i > 0$, and let $\mu \coloneqq \min_{i \in \tom} \mu_i$.
              Then, we get
              \[
                  u_0(x) \le w_\mu(x) \forallcondition{x \in \setR^n}.
              \] \label{enum: u_0 <= u_mu}
        \item For all~$x \in \setR^n$ and~$r \ge \ell > 0$, it follows that
            \[
                  w_r(x) \le w_\ell(x) \le \frac{r}{\ell} w_r(x).
              \]
              \label{enum: u_r w_ell}
    \end{enumerate}
    \end{theorem}

    \section{Multiobjective proximal-PL inequality} \label{sec: PL}
    For~\cref{eq: MOP} with~$m = 1$ and~$g_1 = 0$, we say that the objective function satisfies the Polyak-{\L}ojasiewicz (PL) inequality~\cite{Karimi2016,Polyak1963} if the following holds for some~$\tau > 0$:
    \begin{equation} \label{eq: PL}
        \frac{1}{2} \norm*{ \nabla f_1(x) }^2 \ge \tau (f_1(x) - f_1^\ast) \forallcondition{x \in \setR^n},
    \end{equation}
    where~$f_1^\ast$ is the optimal value~$f_1$.
    Moreover, for~\cref{eq: MOP} the proximal-PL inequality~\cite{Karimi2016} is defined by
\begin{equation} \label{eq: proximal PL}
    \frac{1}{2} \mathcal{D}_{g_1}(x, L_1) \ge \tau (F_1(x) - F_1^\ast) \forallcondition{x \in \setR^n}
\end{equation}
for some~$\tau > 0$, where~$F_1^\ast$ is the optimal value of~$F_1$ and~$\mathcal{D}_{g_1}$ is defined by~\cref{eq: D_g}.
Note that if~$g_1 = 0$, then~\cref{eq: proximal PL} reduces to~\cref{eq: PL}.
These conditions imply that the objective function grows quadratically as the gradient increases and are satisfied when~$f_1$ is strongly convex particularly.
It is known that~\cref{eq: proximal PL} is sufficient for the linear convergence rates of the proximal gradient method.
We propose an extension of~\cref{eq: proximal PL} to multiobjective optimization~\cref{eq: MOP}.
\begin{definition}
    We say that the multiobjective proximal-PL inequality holds for problem~\cref{eq: MOP} when
    \begin{equation} \label{eq: multiobjective proximal-PL}
        L w_L(x) \ge \tau u_0(x) \forallcondition{x \in \setR^n}
    \end{equation}
    for some~$\tau > 0$, where~$L$ is defined by~\cref{eq: L}.
\end{definition}
We note that if~$m = 1$, \cref{eq: multiobjective proximal-PL} reduces to the proximal-PL inequality~\cref{eq: proximal PL}.
We can also generalize the PL inequality~\cref{eq: PL} by setting~$g_i = 0$ in~\cref{eq: multiobjective proximal-PL}.
We state below some sufficient conditions for~\cref{eq: multiobjective proximal-PL}.
\begin{proposition} \label{prop: sufficient}
    The following statements hold.
    \begin{enumerate}
        \item If each~$f_i$ is strongly convex with modulus~$\mu_i$, then~\cref{eq: multiobjective proximal-PL} holds with~$\tau \coloneqq L / \max \{ L / \mu, 1 \}$, where~$\mu = \min_{i \in \tom} \mu_i$. \label{enum: strong convexity implies PL}

        \item Assume that each~$f_i$ is defined by~$f_i(x) \coloneqq h_i(A_i x)$ with a strongly convex function~$h_i$ and a linear transformation~$A_i$, and each~$g_i$ is the indicator function of a polyhedral set~$\mathcal{X}_i$, i.e.,
            \[
                g_i(x) \coloneqq
                \begin{dcases}
                    0, &\quad \text{if } x \in \mathcal{X}_i, \\
                    \infty, &\quad \text{otherwise}.
                \end{dcases}
            \]
            If~$\min_{x \in \setR^n} F_i(x)$ has a nonempty optimal solution set~$X_i^\ast$ for all~$i \in \tom$, then~\cref{eq: multiobjective proximal-PL} holds with some constant~$\tau$. \label{enum: strong convexity composed with linear implies PL}
    \end{enumerate}
\end{proposition}
\begin{proof}
    \ref{enum: strong convexity implies PL}:
    Since~$f_i$ is strongly convex, Theorem~\ref{thm: merit relation}~\ref{enum: u_0 <= u_mu} gives
    \[
        u_0(x) \le w_{\mu}(x) \forallcondition{x \in \setR^n}.
    \]
    Combining the above inequality with Theorem~\ref{thm: merit relation}~\ref{enum: u_r w_ell}, we get
    \[
        u_0(x) \le \max \left\{ \frac{L}{\mu}, 1 \right\} w_L (x) \forallcondition{x \in \setR^n},
    \]
    which means
    \[
        L w_L(x) \ge \frac{L}{ \max \{ L / \mu, 1 \} } u_0(x) \forallcondition{x \in \setR^n}.
    \]

    \ref{enum: strong convexity composed with linear implies PL}:
    Since~$\mathcal{X}_i$ is polyhedral, it can be written as~$\{x \in \setR^n \mid B_i x \le c_i \}$ for a matrix~$B_i$ and a vector~$c_i$.
    Now, we show that for all~$i \in \tom$ there exists some~$z_i$ such that
    \[
    X_i^\ast = \Set*{x \in \setR^n}{B_i x \le c_i \text{ and } A_i x = z_i}
    .\]
    To obtain a contradiction, suppose that there exist~$x^1 \in X_i^\ast$ and~$x^2 \in X_i^\ast$ such that~$A_i x^1 \neq A_i x^2$.
    Clearly, we have~$f_i(x^1) = f_i(x^2)$.
    Since each~$h_i$ is strongly convex, we get
    \[
        \begin{split}
            f_i(x^1) &= 0.5 f_i(x^1) + 0.5 f_i(x^2) = 0.5 h_i(A_i x^1) + 0.5 h_i(A_i x^2) \\
                     &> h_i(A_i(0.5 x^1 + 0.5 x^2)) = f_i(0.5 x^1 + 0.5 x^2)
        ,\end{split}
    \]
    which contradicts the fact that~$x^1 \in X_i^\ast$.
    Therefore, we can use Hoffman's error bound~\cite{Hoffman1952}, and so there exists some~$\rho_i > 0$ such that for any~$x \in \setR^n$ there exists~$x_i^\ast \in X_i^\ast$ with
    \[
        \norm*{x - x_i^\ast} \le \rho_i \norm*{ \max \left\{ \begin{bmatrix} B_i \\ A_i \\ - A_i \end{bmatrix} x - \begin{bmatrix} c_i \\ z_i \\ - z_i \end{bmatrix}, 0 \right\} }
    .\]
    Note that the~$\max$ operator on the right-hand side is taken componentwise.
    Since~$B_i x - c_i \le 0$ for all~$x \in \dom F$, we have
    \[
        \norm*{ x - x_i^\ast } \le \rho_i \norm*{ \max \left\{ \begin{bmatrix} A_i \\ - A_i \end{bmatrix} x - \begin{bmatrix} z_i \\ - z_i \end{bmatrix}, 0 \right\}}  \forallcondition{x \in \dom F},
    \]
    which yields
    \[
        \norm*{ x - x_i^\ast }^2 \le \rho_i^2 \norm*{ A_i x - z_i }^2 \forallcondition{x \in \dom F}.
    \]
    Let~$P_i(x) \coloneqq \argmin_{y \in X_i^\ast} \norm*{x - y}$.
    Since~$P_i(x) \in X_i^\ast$, we get~$z_i = A_i P_i(x)$.
    Then, we have
    \begin{equation} \label{eq: hoffman}
        \norm*{ x - P_i(x) }^2 \le \norm*{ x - x_i^\ast }^2 \le \rho_i^2 \norm*{ A_i (x - P_i(x))}^2 \forallcondition{x \in \dom F}.
    \end{equation}
    Now, suppose that~$x \in \dom F$.
    From the definition~\cref{eq: u_0} of~$u_0$, we have
    \[
        \begin{split}
            u_0(x) &= \sup_{y \in \setR^n} \min_{i \in \tom} \{ F_i(x) - F_i(y) \} \\
                   &\le \min_{i \in \tom} \sup_{y \in \setR^n} \{ F_i(x) - F_i(y) \} = \min_{i \in \tom} \{ F_i(x) - F_i(P_i(x)) \},
        \end{split}
    \]
    where the last equality holds because~$P_i(x) \in X_i^\ast$ and so~$P_i(x) = \argmin_{x \in \setR^n} F_i(x)$.
    Assuming that each~$h_i$ is strongly convex with modulus~$\sigma_i > 0$, it follows that
    \begin{multline*}
        u_0(x) \le \min_{i \in \tom} \Bigl\{\nabla h_i (A_i x)^\T (A_i (x - P_i(x))) \\
        + g_i(x) - g_i(P_i(x)) - \frac{\sigma_i}{2} \norm*{A_i (x - P_i(x))}^2\Bigr\} \\
        = \min_{i \in \tom}\Bigl\{ \nabla f_i(x)^\T (x - P_i(x)) + g_i(x) - g_i(P_i(x)) - \frac{\sigma_i}{2} \norm*{A_i (x - P_i(x))}^2\Bigr\}
    .\end{multline*}
    Applying~\cref{eq: hoffman} to the above inequality, we obtain
    \begin{multline*}
        u_0(x) \\
        \begin{aligned}
        &\le \min_{i \in \tom} \left\{\nabla f_i(x)^\T (x - P_i(x)) + g_i(x) - g_i(P_i(x)) - \frac{\sigma_i}{2 \rho_i^2} \norm*{x - P_i(x)}^2\right\} \\
        &\le \min_{i \in \tom} \left\{\nabla f_i(x)^\T (x - P_i(x)) + g_i(x) - g_i(P_i(x)) - \frac{\upsilon}{2} \norm*{ x - P_i(x)}^2\right\}
        ,\end{aligned}
    \end{multline*}
    where~$\upsilon \coloneqq \max_{i \in \tom} \sigma_i / \rho_i^2$.
    Let~$\mathcal{E} \coloneqq \{e \in \setR^m \mid \sum_{i = 1}^m e_i = 1, e_i \ge 0\}$.
    Since~$\min_{i \in \tom} q_i = \min_{e \in \mathcal{E}} \sum_{i = 1}^m e_i q_i$ for any~$q \in \setR^m$, we get
    \begin{multline*}
        u_0(x) \\
        \begin{aligned}
             &\le \min_{e \in \mathcal{E}} \sum_{i = 1}^m e_i \left\{ \nabla f_i(x)^\T (x - P_i(x)) + g_i(x) - g_i(P_i(x)) - \frac{\upsilon}{2} \norm*{ x - P_i(x)}^2 \right\} \\
                   &\le \min_{e \in \mathcal{E}} \sup_{y \in \setR^n} \sum_{i = 1}^m e_i \left\{ \nabla f_i(x)^\T (x - y) + g_i(x) - g_i(y) - \frac{\upsilon}{2} \norm*{ x - y}^2 \right\} \\
                   &= \sup_{y \in \setR^n} \min_{e \in \mathcal{E}} \sum_{i = 1}^m e_i \left\{ \nabla f_i(x)^\T (x - y) + g_i(x) - g_i(y) - \frac{\upsilon}{2} \norm*{ x - y}^2 \right\} \\
                   &= \sup_{y \in \setR^n} \min_{i \in \tom} \left\{ \nabla f_i(x)^\T (x - y) + g_i(x) - g_i(y) - \frac{\upsilon}{2} \norm*{ x - y}^2 \right\}
               = w_\upsilon (x),
        \end{aligned}
    \end{multline*}
    where the first equality follows from Sion's minimax theorem~\cite{Sion1958}, and the third equality comes from the definition~\cref{eq: w_ell} of~$w_\upsilon$.
    Thus, Theorem~\ref{thm: merit relation}~\ref{enum: u_r w_ell} gives
    \[
        u_0(x) \le \max \left\{\frac{L}{\upsilon}, 1\right\} w_L(x).
    \]
    Multiplying both sides by~$L / \max \set{L / \upsilon, 1}$, the proof is complete.
\end{proof}
In Section~\ref{sec: PL case}, we show that~$\{u_0(x^k)\}$ converges to zero linearly with the multiobjective proximal gradient method under the assumption~\cref{eq: multiobjective proximal-PL}.

\section{Convergence rates of the multiobjective proximal gradient method} \label{sec: convergence rates}
Let us now analyze the convergence rates of the multiobjective proximal gradient method~\cite{Tanabe2019}.
Define~$\psi_x \colon \setR^n \to \setR$ by
\begin{equation} \label{eq: psi}
    \psi_x(d) \coloneqq \max_{i \in \tom} \left\{ \nabla f_i(x)^\T d + g_i(x + d) - g_i(x) \right\}.
\end{equation}
The multiobjective proximal gradient method generates a sequence~$\{ x^k \}$ iteratively with the following procedure:
\[
    x^{k + 1} \coloneqq x^k + d^k,
\]
where~$d^k$ is a search direction.
At every iteration~$k$, we define this~$d^k$ by solving
\begin{equation} \label{eq: subprob}
    d^k \coloneqq \argmin_{d \in \setR^n} \left\{ \psi_{x^k}(d) + \frac{\ell}{2} \norm*{ d }^2 \right\},
\end{equation}
with a positive constant~$\ell > 0$.
Note that we have
\begin{equation} \label{eq: psi w_ell}
    \psi_{x^k}(d^k) + \frac{\ell}{2} \norm*{ d^k }^2 = - w_\ell(x^k),
\end{equation}
where~$w_\ell$ is defined by~\cref{eq: w_ell}.
To keep the paper self-contained, we first recall the algorithm as follows.
\begin{algorithm}[The multiobjective proximal gradient method] \label{alg: proximal_gradient}
	\leavevmode \par
	\begin{enumerate}
		\setlength{\leftskip}{2em}
		\item Choose $\ell > L,\ x^0 \in \setR^n$ and set $k \coloneqq 0$.
        \item Compute $d^k$ by solving the subproblem~\cref{eq: subprob}. \label{enum:subprob_Lipschitz}
		\item If $d^k = 0$, then stop.
		\item Set $x^{k+1} \coloneqq x^k + d^k$, $k \coloneqq k + 1$, and go to \ref{enum:subprob_Lipschitz}. \label{enum:next_x_Lipschitz}
	\end{enumerate}
\end{algorithm}
We suppose that the algorithm generates an infinite sequence of iterates from now on.
The following result shows an important property for~$\psi_x$.
\begin{lemma}~\cite[Lemma~4.1]{Tanabe2019} \label{lem: psi property}
    Let~$\{ d^k \}$ be generated by \cref{alg: proximal_gradient} and recall the definition~\cref{eq: psi} of~$\psi_x$. Then, we have
    \[
        \psi_{x^k}(d^k) \le - \ell \norm*{d^k}^2 \forallcondition{k}.
    \]
\end{lemma}
If~$\ell > L$, from the so-called descent lemma~\cite[Proposition~A.24]{Bertsekas1999}, for all $i \in \tom$ we have
\begin{equation} \label{eq: descent}
    F_i(x^{k + 1}) - F_i(x^k) \le \nabla f_i(x^k)^\T d^k + g_i(x^{k + 1}) - g_i(x^k) + \frac{\ell}{2}\norm*{d^k}.
\end{equation}
The right-hand side of~\cref{eq: descent} is less than zero since~$d^k$ is the optimal solution of~\cref{eq: subprob}, so we get
\begin{equation} \label{eq: nonincreasing}
    F_i(x^{k + 1}) \le F_i(x^k).
\end{equation}
\begin{remark}
    When the Lipschitz constant~$L$ is unknown or incomputable, we can use~$\ell$ for~\cref{eq: subprob} calculated by backtracking instead, i.e., we can set the initial value of~$\ell$ appropriately and multiply~$\ell$ by a prespecified scalar~$\gamma > 1$ at each iteration until~\cref{eq: descent} is satisfied.
    Since~$L$ is finite, the backtracking only requires a finite number of steps.
\end{remark}

\subsection{The non-convex case}
When $m = 1$ and~$F_1$ is not convex for~\cref{eq: MOP},~$\set{\norm{ x^k - \prox_{g_1/L_1} (x^k - \nabla f_1(x^k) / L_1) }}$ converges to zero with rate~$O(\sqrt{1 / k})$ if the proximal gradient method is applied, where the proximal operator~``$\prox$'' is defined by~\cref{eq: D_g}~\cite{Beck2017}.
As in Remark~\ref{rem: merit}, note that when~$m = 1$ we have~$w_{L_1}(x^k) \ge (L_1 / 2) \norm{ x^k - \prox_{g_1/L_1} (x^k - \nabla f_1(x^k) / L_1) }^2$.
Now, in the multiobjective context, we show below that~$\set{ \sqrt{w_1(x^k)} }$ still converges to zero with rate~$O(\sqrt{1 / k})$ with Algorithm~\ref{alg: proximal_gradient}.
\begin{theorem} \label{thm: rate nonconvex}
    Suppose that there exists some nonempty set $\mathcal{J} \subseteq \tom$ such that if $i \in \mathcal{J}$ then $F_i(x)$ has a lower bound~$F_i^{\min}$ for all~$x \in \setR^n$.
    Let $F^{\min} \coloneqq \min_{i \in \mathcal{J}} F_i^{\min}$ and $F_0^{\max} \coloneqq \max_{i \in \tom} F_i(x^0)$.
	Then, \cref{alg: proximal_gradient} generates a sequence~$\{x^k\}$ such that
    \[
		\min_{0 \le j \le k - 1} w_1(x^j) \le \frac{(F_0^{\max} - F^{\min}) \max \set*{ 1, \ell }}{k}
    .\]
\end{theorem}
\begin{proof}
    Let $i \in \mathcal{J}$. From~\cref{eq: descent}, we have
    \begin{multline*}
        F_i(x^{k + 1}) - F_i(x^k) \le \nabla f_i(x^{k})^\T d^k + g_i(x^{k + 1}) - g_i(x^k) + \frac{\ell}{2} \norm*{ d^k}^2 \\
        \le \max_{i \in \tom} \left\{ \nabla f_i(x^{k})^\T d^k + g_i(x^{k + 1}) - g_i(x^k) + \frac{\ell}{2} \norm*{ d^k}^2 \right\}
        = - w_\ell(x^k)
    ,\end{multline*}
    where the equality follows from~\cref{eq: psi} and~\cref{eq: psi w_ell}.
    Adding up the above inequality from~$k = 0$ to~$k = \tilde{k} - 1$ yields that
    \[
        F_i(x^{\tilde{k}}) - F_i(x^0) \le - \sum_{k = 0}^{\tilde{k} - 1} w_\ell(x^k)
		\le - \tilde{k} \min_{0 \le k \le \tilde{k} - 1} w_\ell(x^k).
    \]
    From the definitions of~$F^{\min}$ and~$F_0^{\max}$, we obtain
    \[
        \min_{0 \le k \le \tilde{k} - 1} w_\ell(x^k) \le \frac{F_0^{\max} - F^{\min}}{\tilde{k}}.
    \]
    Finally, from Theorem~\ref{thm: merit relation}~\ref{enum: u_r w_ell}, we get
	\[
        \min_{0 \le k \le \tilde{k} - 1} w_1(x^k) \le \frac{(F_0^{\max} - F^{\min}) \max \{ 1, \ell \}}{\tilde{k}}.
	\]
\end{proof}
\begin{remark}
    When~$g_i = 0$ for all~$i$, references~\cite{Calderon2020,Fliege2019,Grapiglia2015} present the convergence rate of various multiobjective optimization methods.
    However, as we have mentioned in the introduction, they all evaluate the convergence rate with measures that depend on the subproblems or variables used in their algorithms.
    This means that the comparison in terms of complexity between different methods is not easy by using those measures.
    However, \cref{thm: rate nonconvex} analyzes the convergence rate using the merit function~$w_1$, which can be defined uniformly by~\cref{eq: w_ell} for multiobjective optimization problems with a structure like~\cref{eq: f + g}.
\end{remark}

\subsection{The convex case}
For \cref{eq: MOP} with~$m = 1$, if~$f_1$ is convex, then $\{ F_1(x^k) - F_1^\ast \}$ converges to zero with rate~$O(1/k)$ using the proximal gradient method, where~$F_1^\ast$ is the optimal objective value of~$F_1$~\cite{Beck2009}.
In this subsection, we show how fast~$\{u_0(x^k)\}$ converges to zero with Algorithm~\ref{alg: proximal_gradient}.
Let us start by proving the following lemma.
Note, however, that we state it with $f_i$ and $g_i$ having general (nonnegative) convexity parameters\footnote{We say that~$h \colon \setR^n \to \setR \cup \set*{\infty}$ has a \emph{convexity parameter}~$\varsigma \in \setR$ if~$h(\alpha x + (1 - \alpha)y) \le \alpha h(x) + (1 - \alpha)h(y) - (1 / 2)\alpha (1 - \alpha) \varsigma \norm*{x - y}^2$ holds for all~$x, y \in \setR^n$ and~$\alpha \in [0, 1]$. When~$\varsigma > 0$, we call~$h$ strongly convex. Note that this definition allows non-convex cases, i.e.,~$\varsigma < 0$ can also be considered.}, but they turn out to be zero in this subsection.
\begin{lemma} \label{lem:without_line_lem}
	Let $f_i$ and $g_i$ have convexity parameters~$\mu_i \in \setR$ and $\nu_i \in \setR$, respectively, and write $\mu \coloneqq \min_{i \in \tom} \mu_i$ and $\nu \coloneqq \min_{i \in \tom} \nu_i$.
    Then, for all $x \in \setR^n$ it follows that
    \begin{multline*}
        \sum_{i = 1}^m \lambda_i^k \left(F_i(x^{k + 1}) - F_i(x)\right) \le \frac{\ell}{2} \left( \norm*{ x^k - x}^2 - \norm*{ x^{k + 1} - x}^2 \right) \\
        - \frac{\nu}{2} \norm*{ x^{k + 1} - x }^2 - \frac{\mu}{2} \norm*{ x^k - x}^2,
    \end{multline*}
	where $\lambda_i^k$ satisfies the following conditions:
    \begin{enumerate}
        \item There exists $\eta_i^k \in \partial g_i(x^k + d^k)$ such that $\sum_{i = 1}^m \lambda_i^k (\nabla f_i(x^k) + \eta_i^k) + \ell d^k = 0$,
        \item $\sum_{i = 1}^m \lambda_i^k = 1, \lambda_i^k \ge 0 \left(i \in \mathcal{I}_{x^k}(d^k)\right)$ and $\lambda_i^k = 0 \left(i \notin \mathcal{I}_{x^k}(d^k)\right)$,
    \end{enumerate}
	where~$\mathcal{I}_x(d) \coloneqq \Set*{i \in \tom}{\psi_x(d) = \nabla f_i(x)^\T d + g_i(x + d) - g_i(x)}$.
\end{lemma}
\begin{proof}
    From~\cref{eq: descent}, we have
	\[
        F_i(x^{k + 1}) - F_i(x^k) \le \nabla f_i(x^k)^\T (x^{k + 1} - x^k) + g_i(x^{k + 1}) - g_i(x^k) + \frac{\ell}{2} \norm*{ d^k}^2
	.\]
    The above inequality and convexity of $f_i$ with modulus~$\mu_i$ give
    \[
        \begin{split}
            \MoveEqLeft[.5] F_i(x^{k + 1}) - F_i(x) = ( F_i(x^k) - F_i(x) ) + ( F_i(x^{k + 1}) - F_i(x^k) ) \\
        &\le \left(\nabla f_i(x^k)^\T (x^k - x) - \frac{\mu_i}{2} \norm*{ x^k - x}^2 + g_i(x^k) - g_i(x)\right) \\
        &\quad + \left(\nabla f_i(x^k)^\T (x^{k + 1} - x^k) + g_i(x^{k + 1}) - g_i(x^k) + \frac{\ell}{2} \norm*{ x^{k + 1} - x^k}^2\right) \\
        &\le \nabla f_i(x^k)^\T (x^k + d^k - x) + g_i(x^k + d^k) - g_i(x)  - \frac{\mu}{2} \norm*{ x^k - x}^2 + \frac{\ell}{2} \norm*{ d^k}^2 \\
        &\le (\nabla f_i(x^k) + \eta_i^k)^\T (x^k + d^k - x)
        - \frac{\mu}{2} \norm*{ x^k - x}^2 - \frac{\nu}{2} \norm*{ x^{k + 1} - x }^2 + \frac{\ell}{2} \norm*{d^k}^2
        ,\end{split}
    \]
	where the second inequality follows from the definition of $\mu$ and the fact that $x^{k + 1} = x^k + d^k$, and the last one comes from the convexity of~$g_i$. Multiplying the above inequality by $\lambda_i^k$ and summing for all~$i \in \tom$, the conditions~(i) and (ii) give
    \[
        \begin{split}
            \MoveEqLeft[.5] \sum_{i = 1}^m \lambda_i^k \left( F_i(x^{k + 1}) - F_i(x) \right) \\
        &= {- \ell (d^k)^\T (x^k + d^k - x) - \frac{\mu}{2} \norm*{ x^k - x}^2 - \frac{\nu}{2} \norm*{ x^{k + 1} - x }^2 + \frac{\ell}{2} \norm*{ d^k}^2} \\
        &=  - \frac{\ell}{2} \left( 2 (d^k)^\T (x^k - x) + \norm*{ d^k }^2 \right) - \frac{\mu}{2} \norm*{ x^k - x}^2 - \frac{\nu}{2} \norm*{ x^{k + 1} - x }^2 \\
        &=  \frac{\ell}{2} \left( \norm*{ x^k - x}^2 - \norm*{ x^{k + 1} - x}^2 \right) - \frac{\mu}{2} \norm*{ x^k - x}^2 - \frac{\nu}{2} \norm*{ x^{k + 1} - x }^2
        .\end{split}
    \]
\end{proof}

Now, we show that~$\set{u_0(x^k)}$ converges to zero with rate~$O(1 / k)$ with Algorithm~\ref{alg: proximal_gradient} under the following assumption.
\begin{assumption} \label{asm: boundedness}
    Let~$X^\ast$ be the set of weakly Pareto optimal points for~\cref{eq: MOP}, and define the level set of~$F$ for~$\alpha \in \setR^m$ by~$\Omega_F(\alpha) \coloneqq \{ x \in S \mid F(x) \le \alpha \}$.
    Then, for all~$x \in \Omega_F(F(x^0))$ there exists~$x^\ast \in X^\ast$ such that~$F(x^\ast) \le F(x)$ and
    \[
        R \coloneqq \sup_{F^\ast \in F(X^\ast \cap \Omega_F(F(x^0)))} \inf_{x \in F^{-1}(\set*{F^\ast})} \norm*{x - x^0}^2 < \infty
    .\]
\end{assumption}
\begin{remark}
    \begin{enumerate}[wide, labelwidth=!, labelindent=0pt]
        \item In single-objective cases, Assumption~\ref{asm: boundedness} is valid if the optimization problem has at least one optimal solution; when~$m = 1$,~$X^\ast$ coincides with the optimal solution set, and the equality~$X^\ast \cap \Omega_F(F(x^0)) = X^\ast$ holds, so we have~$R = \inf_{x \in X^\ast} \norm*{x - x^0}^2 < \infty$.
        \item When the level set~$\Omega_F(F(x^0))$ is bounded, Assumption~\ref{asm: boundedness} is also satisfied.
            For example, this is the case when~$F_i$ is strongly convex for at least one~$i$.
    \end{enumerate}
\end{remark}
\begin{theorem}
	Assume that~$F_i$ is convex for all~$i \in \tom$.
    Under Assumption~\ref{asm: boundedness}, Algorithm~\ref{alg: proximal_gradient} generates a sequence~$\{ x^k \}$ such that
    \[
        u_0(x^k) \le \frac{\ell R}{2 k} \forallcondition{k \ge 1}.
    \]
\end{theorem}
\begin{proof}
	From Lemma~\ref{lem:without_line_lem} and the convexity of~$f_i$ and~$g_i$, for all~$x \in \setR^n$ we have
    \[
		\sum_{i = 1}^m \lambda_i^k \left(F_i(x^{k + 1}) - F_i(x)\right) \le \frac{\ell}{2} \left( \norm*{x^k - x}^2 - \norm*{x^{k + 1} - x}^2 \right)
    .\]
    Adding up the above inequality from~$k = 0$ to~$k = \hat{k}$, we obtain
    \[
        \begin{split}
            \sum_{k = 0}^{\hat{k}} \sum_{i = 1}^m \lambda_i^k \left(F_i(x^{k + 1}) - F_i(x)\right) &\le \frac{\ell}{2} \left( \norm*{ x^0 - x}^2 - \norm*{ x^{\hat{k} + 1} - x}^2 \right) \\
        &\le \frac{\ell}{2} \norm*{x^0 - x}^2
        .\end{split}
    \]
    Since~$F_i(x^{\hat{k} + 1}) \le F_i(x^{k + 1})$ for all~$k \le \hat{k}$ (see~\cref{eq: nonincreasing}), we get
    \[
        \sum_{k = 0}^{\hat{k}} \sum_{i = 1}^m \lambda_i^k \left(F_i(x^{\hat{k} + 1}) - F_i(x)\right) \le \frac{\ell}{2} \norm*{x^0 - x}^2.
    \]
    Let~$\bar{\lambda}_i^{\hat{k}} \coloneqq \sum_{k = 0}^{\hat{k}} \lambda_i^k / (\hat{k} + 1)$. Then, it follows that
	\[
        \sum_{i = 1}^m \bar{\lambda}_i^{\hat{k}} \left(F_i(x^{\hat{k} + 1}) - F_i(x)\right)  \le \frac{\ell}{2 (\hat{k} + 1)} \norm*{x^0 - x}^2.
	\]
    Since~$\bar{\lambda}_i^{\hat{k}} \ge 0$ and $\sum_{i = 1}^m \bar{\lambda}_i^{\hat{k}} = 1$, we see that
	\[
        \min_{i \in \tom} \left( F_i(x^{\hat{k} + 1}) - F_i(x) \right) \le \frac{\ell}{2 (\hat{k} + 1)} \norm*{x^0 - x}^2.
	\]
    Therefore, we get
	\[
        \sup_{F^\ast \in F(X^\ast \cap \Omega_F(F(x^0)))} \inf_{x \in F^{-1}(\{F^\ast\})} \min_{i \in \tom} \left( F_i(x^{\hat{k} + 1}) - F_i(x) \right) \le \frac{\ell R}{2 (\hat{k} + 1)}.
	\]
    Thus, we obtain
    \[
        \sup_{F^\ast \in F(X^\ast \cap \Omega_F(F(x^0)))} \min_{i \in \tom} \left( F_i(x^{\hat{k} + 1}) - F_i^\ast \right) \le \frac{\ell R}{2 (\hat{k} + 1)},
    \]
    which gives
    \begin{equation} \label{eq: pareto level bound}
        \sup_{x \in X^\ast \cap \Omega_F(F(x^0))} \min_{i \in \tom} \left( F_i(x^{\hat{k} + 1}) - F_i(x) \right) \le \frac{\ell R}{2 (\hat{k} + 1)}.
    \end{equation}
    Now, the inequality~$F_i(x^k) \le F_i(x^0)$ from~\cref{eq: nonincreasing} gives
    \[
        \begin{split}
            \MoveEqLeft[.5] \sup_{x \in \Omega_F(F(x^0))} \min_{i \in \tom} \left( F_i(x^{\hat{k} + 1}) - F_i(x) \right) \\
        &= \sup_{x \in \Omega_F(F(x^{\hat{k} + 1}))} \min_{i \in \tom} \left( F_i(x^{\hat{k} + 1}) - F_i(x) \right)
        ,\end{split}
    \]
    so we have
    \begin{equation} \label{eq: level merit}
        \begin{split}
            \MoveEqLeft \sup_{x \in \Omega_F(F(x^0))} \min_{i \in \tom} \left( F_i(x^{\hat{k} + 1}) - F_i(x) \right) \\
        &= \sup_{x \in \setR^n} \min_{i \in \tom} \left( F_i(x^{\hat{k} + 1}) - F_i(x) \right).
        \end{split}
    \end{equation}
    Moreover, from the assumption that for all~$x \in \Omega_F(F(x^0))$ there exists~$x^\ast \in X^\ast$ such that~$F(x^\ast) \le F(x)$, it follows that
    \begin{equation} \label{eq: pareto level}
        \begin{split}
            \MoveEqLeft \sup_{x \in X^\ast \cap \Omega_F(F(x^0))} \min_{i \in \tom} \left( F_i(x^{\hat{k} + 1}) - F_i(x) \right) \\
            &= \sup_{x \in \Omega_F(F(x^0))} \min_{i \in \tom} \left( F_i(x^{\hat{k} + 1}) - F_i(x) \right)
        .\end{split}
    \end{equation}
    Finally, from~\cref{eq: pareto level bound,eq: level merit,eq: pareto level} we conclude that
    \[
        u_0(x^{\hat{k} + 1}) \coloneqq \sup_{x \in \setR^n} \min_{i \in \tom} \left( F_i(x^{\hat{k} + 1}) - F_i(x) \right) \le \frac{\ell R}{2 (\hat{k} + 1)}.
    \]
\end{proof}

\subsection{The strongly convex case}
For~\cref{eq: MOP} with~$m = 1$, it is known that~$\{x^k\}$ converges linearly to the optimal point when~$f_1$ is strongly convex~\cite{Beck2017}.
Now, we show that the same result holds for Algorithm~\ref{alg: proximal_gradient}.
\begin{theorem} \label{thm: strongly convex}
    Let $f_i$ and $g_i$ have convexity parameters~$\mu_i \in \setR$ and $\nu_i \in \setR$, respectively, and write $\mu \coloneqq \min_{i \in \tom} \mu_i$ and $\nu \coloneqq \min_{i \in \tom} \nu_i$. If $\ell > L$, then there exists a Pareto optimal point~$x^\ast \in \setR^n$ such that for each iteration~$k$,
	\[
		\norm*{ x^{k + 1} - x^\ast } \le \sqrt{\frac{\ell - \mu}{\ell + \nu}} \norm*{x^k - x^\ast}.
	\]
	Thus, we have
	\[
		\norm*{ x^k - x^\ast } \le \left( \sqrt{\frac{\ell - \mu}{\ell + \nu}} \right)^k \norm*{x^0 - x^\ast}.
	\]
\end{theorem}
\begin{proof}
    Since each~$F_i$ is strongly convex, the level set of every~$F_i$ is bounded.
    Thus, $\{ x^k \}$ has an accumulation point~$x^\ast \in \setR^n$.
    Note that~$x^\ast$ is a Pareto optimal point~\cite[Lemma~2.2~and~Theorem~4.3]{Tanabe2019}.
    Now, from Lemma~\ref{lem:without_line_lem}, we have
    \[
        \begin{split}
            \sum_{i = 1}^m \lambda_i^k \left(F_i(x^{k + 1}) - F_i(x^\ast)\right) \le{}& \frac{\ell}{2} \left( \norm*{ x^k - x^\ast }^2 - \norm*{ x^{k + 1} - x^\ast }^2 \right) \\
    &- \frac{\mu}{2} \norm*{ x^k - x^\ast }^2 - \frac{\nu}{2} \norm*{ x^{k + 1} - x^\ast }^2.
        \end{split}
    \]
    Since the left-hand side is nonnegative because of~\cref{eq: psi},~\cref{eq: descent}, and \cref{lem: psi property}, we obtain
	\[
		0 \le \frac{\ell}{2} \left( \norm*{x^k - x^\ast }^2 - \norm*{x^{k + 1} - x^\ast }^2 \right) - \frac{\mu}{2} \norm*{x^k - x^\ast }^2 - \frac{\nu}{2} \norm*{ x^{k + 1} - x^\ast }^2,
	\]
	which is equivalent to
	\[
		\norm*{ x^{k + 1} - x^\ast } \le \sqrt{\frac{\ell - \mu}{\ell + \nu}} \norm*{x^k - x^\ast}.
	\]
\end{proof}

\subsection{The case that the multiobjective proximal-PL inequality is assumed}
\label{sec: PL case}
For~\cref{eq: MOP} with~$m = 1$ satisfying proximal-PL inequality~\cref{eq: proximal PL}, it is known that~$\{ F_1(x^k) - F_1^\ast \}$ converges linearly to zero with the proximal gradient method~\cite{Karimi2016}.
Now, we show that~$\{u_0(x^k)\}$ converges linearly to zero with \cref{alg: proximal_gradient} for the multiobjective problem~\cref{eq: MOP} that satisfies the multiobjective proximal-PL inequality~\cref{eq: multiobjective proximal-PL}.
\begin{theorem} \label{thm: PL}
    Suppose that~\cref{eq: multiobjective proximal-PL} holds with a constant~$\tau > 0$.
    Then, Algorithm~\ref{alg: proximal_gradient} generates a sequence~$\{ x^k \}$ such that
    \[
        u_0 (x^{k + 1}) \le \left( 1 - \frac{\tau}{\ell} \right) u_0(x^k).
    \]
\end{theorem}
\begin{proof}
    Since~$\nabla f_i$ is Lipschitz continuous with constant~$\ell > L$, for all~$i \in \tom$ we get
    \[
        \begin{split}
            F_i(x^{k + 1}) - F_i(x^k) &\le \nabla f_i(x^k)^\T d^k + g_i(x^{k + 1}) - g_i(x^k) + \frac{\ell}{2}\norm*{ d^k}^2 \\
                                      &\le \max_{i \in \tom} \left\{ \nabla f_i(x^k)^\T d^k + g_i(x^{k + 1}) - g_i(x^k) + \frac{\ell}{2}\norm*{ d^k}^2 \right\} \\
                                      &= - w_\ell(x^k)
        \le - \frac{\tau}{\ell} u_0(x^k)
        ,\end{split}
    \]
    where the equality follows from~\cref{eq: psi w_ell}, and the last inequality comes from~\cref{eq: multiobjective proximal-PL}.
    Then, for all~$x \in \setR^n$ we obtain
    \[
        F_i(x^{k + 1}) - F_i(x) \le F_i(x^k) - F_i(x) - \frac{\tau}{\ell} u_0(x^k)
    .\]
    Therefore, it follows that
    \[
        \sup_{x \in \setR^n} \min_{i \in \tom} \{ F_i(x^{k + 1}) - F_i(x) \} \le \sup_{x \in \setR^n} \min_{i \in \tom} \{ F_i(x^k) - F_i(x) \} - \frac{\tau}{\ell} u_0(x^k),
    \]
    which is equivalent to
    \[
        u_0 (x^{k + 1}) \le \left( 1 - \frac{\tau}{\ell} \right) u_0(x^k).
    \]
\end{proof}

\begin{remark}
    While Theorem~\ref{thm: strongly convex} implies linear convergence of~$\{x^k\}$, we show that~$\{u_0(x^k)\}$ converges to zero linearly in Theorem~\ref{thm: PL}.
    However, from Theorem~\ref{thm: PL} and \cref{prop: sufficient}, we can show that if each~$f_i$ is strongly convex with modulus~$\mu_i > 0$, then it follows that
    \[
        u_0 (x^{k + 1}) \le \left( 1 - \frac{L}{\ell \max \{L / \mu, 1\}} \right) u_0(x^k),
    \]
    with~$\mu \coloneqq \max_{i \in \tom} \mu_i$.
\end{remark}

\section{Conclusion} \label{sec: conclusion}
We presented global convergence rates for the multiobjective proximal gradient method~\cite{Tanabe2019}, matching what we know in scalar optimization for non-convex, convex, and strongly convex problems using merit functions to measure the complexity.
We also proposed the multiobjective proximal-Polyak-{\L}ojasiewicz (proximal-PL) inequality that extended proximal-PL inequality for scalar optimization and established the rate for multiobjective problems that satisfy such inequalities.
The results obtained here can be applied to the multiobjective steepest descent method~\cite{Fliege2000}, too, since the multiobjective proximal gradient method generalizes them.

\section*{Acknowledgements}
    This work was supported by the Grant-in-Aid for Scientific Research (C) (17K00032 and 19K11840) and Grant-in-Aid for JSPS Fellows (20J21961) from the Japan Society for the Promotion of Science.
    We are also grateful to the anonymous referees for their useful comments.

    This version of the article has been accepted for publication, after peer review (when applicable) and is subject to Springer Nature’s AM terms of use, but is not the Version of Record and does not reflect post-acceptance improvements, or any corrections. The Version of Record is available online at: \url{http://dx.doi.org/10.1007/s11590-022-01877-7}

\bibliographystyle{jorsj}
\bibliography{library}

\end{document}